\documentclass[a4paper,12pt,notitlepage]{article}
\usepackage[english]{babel}
\usepackage[utf8]{inputenc}
\usepackage{graphicx}
\usepackage{amsmath}
\usepackage{amssymb}
\usepackage{indentfirst}
\usepackage{color,graphicx,caption}
\usepackage{float}
\usepackage{wrapfig}
\usepackage{enumitem}
\usepackage{cancel}
\usepackage{tikz}
\usepackage{caption}
\usepackage{titling}

\pretitle{\begin{center}\bfseries}
\posttitle{\par\end{center}\vskip 0.5em}
\preauthor{\begin{center}\ttfamily}
\postauthor{\end{center}}
\predate{\par\centering}
\postdate{\par}

\newcommand\blfootnote[1]{%
  \begingroup
  \renewcommand\thefootnote{}\footnote{#1}%
  \addtocounter{footnote}{-1}%
  \endgroup
}

\newtheorem{thm}{Theorem}[section]
\newtheorem{definicao}[thm]{\bf Definition}
\newtheorem{prop}[thm]{\bf Proposition}
\newtheorem{teorema}[thm]{\bf Theorem}
\newtheorem{lema}[thm]{\bf Lemma}
\newtheorem{corolario}[thm]{\bf Corollary}
\newtheorem{observacao}[thm]{\bf Remark}

\newenvironment{prova}[1][\bf Proof: \ ]{\textbf{#1}}{\hfill\rule{2mm}{2mm}}
\newcommand{\R}{\mathbb{R}} %reais
\newcommand{\rn}{\mathbb{R}^n} %r^n
\renewcommand{\d}{\textrm{d}}% d retinho de derivada
\newcommand{\norm}[1]{\left\lVert#1\right\rVert} %norma
\newcommand\supp{\mathrm{supp\,}}

\newcommand{\IRn}[1]{\int_{\R^n} #1\,\,\d x }

\newcommand{\Iab}[3]{\int_{#1}^{#2} #3\,\,\d s }
\newcommand{\Iabt}[3]{\int_{#1}^{#2} #3\,\,\d t }
\newcommand{\Iot}[1]{\int_0^t #1\,\,\d s }

\newcommand{\Iabtau}[3]{\int_{#1}^{#2} #3\,\d \tau } %Int de a a b, dt
\newcommand{\IBR}[1]{\int_{B_R} #1\,\d x }
\newcommand{\ddt}{\dfrac{\mathrm{d}}{\mathrm{d}t}}

\newcommand{\implica}{\Rightarrow}

\title{NON-EXISTENCE OF SOLUTIONS FOR A FOURTH ORDER  STRUCTURALLY DAMPED  EQUATION UNDER NONLINEAR MEMORY EFFECT. }
\author{Luis Gustavo Longen \thanks{Universidade Federal de Santa Catarina, Brazil}}
\date{September, 2024}

\begin{document}

\maketitle

\begin{abstract}
In this paper, we study a mixed wave-plate equation with rotational inertia, fractional damping and memory non-linearity. This research is a non-existence counterpart to a paper by D'Abbicco and Longen \cite{longen}, in search of the critical exponent for the global in-time existence of small data solutions to the Cauchy problem for the equation
\[u_{tt}-\Delta u_{tt}+\Delta^2 u-\Delta u+(-\Delta)^\theta u_t=\int_0^t(t-s)^{-\gamma}|u(s,x)|^p\,\d s,\qquad t\in \R_+^\ast, x\in\rn,\]
with $\theta\in[0,\frac{1}{2}),\gamma\in(0,1)$ and $p>1.$ We employ a modified test function argument to show that there are arbitrarily small initial data for which there are no solutions for the problem in the subcritical case $p<\max\lbrace \gamma^{-1},p_c\rbrace,$ where
\[p_c(n,\gamma,\theta):=1+\frac{2\left(1+(1-\gamma)(1-\theta)\right)}{\left(n-2+2\gamma(1-\theta)\right)_+}.\]
Additionally, we show non-existence of global solutions in the critical case $p=\max\lbrace \gamma^{-1},p_c\rbrace$ when $\gamma>\frac{n-2}{n},$ that is, when $p_c>\gamma^{-1}.$

\end{abstract}

\blfootnote{2020 Mathematics Subject Classification 35B33, 35L71, 74K20.}
\blfootnote{ \emph{Key words and phrases:} nonlinear evolution equations, critical exponent, global-in-time solutions, fractional Laplacian, test function method}
\blfootnote{\textbf{Correspondence:} Luis Gustavo Longen, Departamento de Matemática, Universidade Federal de Santa Catarina, R. Eng. Agronômico Andrei Cristian Ferreira, s/n - Trindade, Florianópolis, Santa Catarina, Brazil}
\blfootnote{Email: luisglongen@gmail.com}

%begin{center}
%\textbf{NON-EXISTENCE OF SOLUTIONS FOR A PLATE-LIKE EQUATION WITH FRACTIONAL DAMPING AND MEMORY NON-LINEARITY}

%\vspace{20pt}
%LUIS GUSTAVO LONGEN
%\end{center}

%ABSTRACT. lalala

\section{Introduction}

In this work, we continue an investigation of the critical exponent for  global solutions for a class of semilinear evolution equations with fractional damping and a nonlinear memory term. More precisely, we consider the following Cauchy problem

\begin{equation}\begin{cases}
 u_{tt}-\Delta u_{tt}+\Delta^2 u-\Delta u+(-\Delta)^\theta u_t=F(t,u),\,\,\,\,\,t\ge0,\,\,x\in\rn.\\
(u,u_t)(0,x)=(u_0,u_1)(x),
\end{cases} \tag{P} \label{P}\end{equation}
with $\theta\in[0,\frac{1}{2}),\gamma\in(0,1)$ and $p>1,$ where the operator
\[F(t,u)=\int_0^t(t-s)^{-\gamma}|u(s,x)|^p\,\d s\]
is referred to as the nonlinear memory of the function $u$.

For the whole text, the fractional Laplace operator, denoted by $(-\Delta)^\theta$ in the damping term of $\eqref{P}$ is to be understood as the unique extension to the Schwartz space $\mathcal{S}$ of tempered distributions, of the operator
\[(-\Delta)^\theta f=\mathcal{F}^{-1}\left(|\xi|^{2\theta} \hat{f}\right),\]
where $\hat{f}=\mathcal{F}f$ denotes the Fourier transform of $f,$
\[\mathcal{F}f(\xi)=\hat{f}(\xi)=(2\pi)^{-\frac n2} \IRn{e^{-ix\cdot \xi}|f(x)|},\qquad \xi\in\rn.\]

For \emph{critical exponent} we understand a threshold value $\overline{p}:=\overline{p}(n,\gamma,\theta),$ $\overline{p}>1,$ such that
\begin{itemize}
    \item If $p>\overline{p},$ then there exist Global in-time Small Data Solutions to \eqref{P}, for a suitable choice of initial data and adequate solution spaces.
    \item If $1<p<\overline{p},$ then there can be found arbitrarily small initial data such that there is no global in-time solution to \eqref{P}.
\end{itemize}

We observe that $(-\Delta)^{\theta}$ is nonnegative and self-adjoint in $L^2(\rn).$ Therefore, its square root, $(-\Delta)^{\frac{\theta}{2}}$, is well-defined, and also self-adjoint. If one considers the linear associated problem to \eqref{P}, that is,
\begin{equation}\begin{cases}
 u_{tt}-\Delta u_{tt}+\Delta^2 u-\Delta u+(-\Delta)^\theta u_t=0,\,\,\,\,\,t\ge0,\,\,x\in\rn.\\
(u,u_t)(0,x)=(u_0,u_1)(x),
\end{cases} \tag{LP} \label{PL}\end{equation}
by formally taking the inner product in $L^2(\rn)$ with $u_t,$ it becomes clear that
\begin{equation}
    \ddt E(t) + 2\norm{(-\Delta)^\frac{\theta}{2}u_t}^2_{L^2}=0, \label{eq_energia}
\end{equation}
where the associated energy to the problem $E(t)$ is
\begin{equation}
    E(t):=\left(\norm{ u_t}^2_{L^2}+\norm{\nabla  u_t}^2_{L^2}+\norm{\Delta u}^2_{L^2}\norm{\nabla  u}^2_{L^2} \right).
\end{equation}
Therefore, the natural space for solutions to problems $\eqref{P}$ and $\eqref{PL}$ is
\[u\in \mathcal{C}\left([0,\infty),H^2(\rn)\right)\cap \mathcal{C}^1\left([0,\infty),H^1(\rn)\right). \]

\subsection{Notation}
We list some notation used in this paper:
\begin{itemize}
    \item the expression $(a)_+$ denotes the positive part of $a,$ that is, $(a)_+=\max\lbrace a,0 \rbrace;$
    \item the expression $f(t)\lesssim g(t)$ denotes that there exist a constant $C>0$ such that $f(t)\le C g(t),$ uniformly with respect to $t$;
    \item for any quantity $x\in \rn,$ we denote  $\langle x \rangle = (1+|x|^2)^{\frac{1}{2}};$
    \item the expression $\lfloor x \rfloor$ denotes the floor function of $x\in\R,$ that is, the greatest integer less than or equal to $x$;

\end{itemize}

\subsection*{Historical Review}

We start reviewing some known historical results concerning the critical exponent, and existence or non-existence of solutions for the supercritical and subcritical cases, respectively.

In his pioneer work in 1966, H. Fujita \cite{fujita} proved that the critical exponent for the semilinear heat equation
\begin{equation}
    u_t-\Delta u=u^p
    \label{eq:fujita}
\end{equation}
is $p_F:=1+\frac2n.$ This is known as the \emph{Fujita exponent.} In 2001, Todorova and Yordanov \cite{todorova} proved that the critical exponent is still $p_F$ for the semilinear heat equation with nonlinearity $F(u)=|u|^p,$ that is,
\begin{equation}
    u_t-\Delta u = |u|^p.  \tag{NLH}
    \label{eq:todorova}
\end{equation}
In 2003, Zhang \cite{zhang} proved the nonexistence for the critical case for \eqref{eq:todorova}.

In 1968, Strauss \cite{strauss1} investigated the semilinear wave equation in 3 dimensions, obtaining another type of critical exponent. This critical exponent was later generalized for higher dimensions by Glassey \cite{glassey1}, \cite{glassey2} . Glassey showed that the critical exponent for the problem
\begin{equation}
    u_{tt}-\Delta u = |u|^p \tag{NLW}
\end{equation}
is $p_S,$ the larger root of the quadratic equation
\[(n-1)p^2-(n+1)p-2=0.\]
This is now known as the \emph{Strauss exponent.} In the following years, several mathematicians have studied semilinear problems of the form
\begin{equation}
    \dfrac{\mathrm{d}^j u}{\mathrm{d}t^j}=Hu+Pu, \qquad j=1,2
\end{equation}
where $H$ is a differential operator acting in the spatial coordinates and $P=P(u,p)$ is a multiplication by a ``potential'' function depending on the unknown function $u$ and the exponent $p$. This much more general equation can be reduced to several interesting models in physics (see \cite{strauss2}), for example the nonlinear Klein-Gordon equation
\begin{equation}
    u_{tt}-\Delta u +m^2u+h(u)=0, \qquad m>0,\,\, x\in\rn, \tag{NLKG}
\end{equation}
where $h$ satisfies the following condition
\begin{equation}
|h'(s)|\lesssim |s|^{p-1}.    \label{eq:cond_h}
\end{equation}
 
The same model can be applied to the nonlinear Schrödinger Equation 
\begin{equation}
    i u_t-\Delta u +h(|u|)\arg(u), \tag{NLS}
\end{equation}
where $h$ satisfies \eqref{eq:cond_h}. Another important application of the same model is the Generalized Korteweg-deVries equation
\begin{equation}
    u_t+u_{xxx}+h(u)_x=0,\qquad x\in\rn. \tag{GKdV}
\end{equation}
Each one of these models have been studied thoroughly, and obtained results have been critical exponents that are either generalizations or shifts of $p_F$ or $p_S.$ Only in 2014, Cazenave, Dickstein and Weissler \cite{cazenave} encountered a completely different kind of critical exponent. Investigating the heat equation with a memory term for nonlinearity,
\begin{equation}
    u_t-\Delta u = \Iot{(t-s)^\gamma |u|^{p-1}u(s)},
\end{equation}
the critical exponent found was $p_\star=\max\lbrace p_\gamma,\gamma^{-1}\rbrace,$ where
\[p_\gamma=1+\dfrac{4-2\gamma}{(n-2+2\gamma)_+}\]
is a Fujita-type exponent. The same kind of result, a competition between a Fujita-type exponent and $\gamma^{-1},$ was obtained for the damped wave equation with the same nonlinearity by D'Abbicco \cite{DA14NODEA} in 2014.

Turning attention to the non-existence of solutions for the subcritical case, the test function method was introduced by Mitidieri and Pohozaev \cite{Mitidieri1}, \cite{Mitidieri2}, \cite{Mitidieri3} in 1998. After this, the method was employed in several works. A small (and far from exhaustive) list of such results follows: D'Ambrosio-Reissig \cite{lucente}, Fujiwara-Ozawa \cite{fuji-ozawa}, Ikeda-Inui \cite{ikeda}, 
Lai-Zhou \cite{lai} and Lin-Nishihara \cite{lin}.

Finally, it is worth noticing that the critical case is very sensitive, in some sense. In fact, even for the classical heat equation, nonexistence of solutions was proved for the critical case only two years after the critical exponent was found, and the techniques applied were quite different from the ones used in the subcritical and supercritical cases. Also, although many works concerning the critical case conclude the non-existence of solutions for $p=p_c,$ it is worth noticing that Ebert, Da Luz and Palma \cite{mairaebertcluz} showed in their work in 2020 that there exist Global Solutions for the critical case, if the initial data has more regularity.

The study of fourth-order evolution partial differential equations is very important for modelling problems in solid mechanics, such as  large deflections in thin plates. In this sense, it is known (as one can see in \cite{lagnese}) that the \emph{full von Karmán dynamical model}, namely
\[\begin{cases}
    U_{tt}=\textrm{Div}\left(\mathcal{C}\left[e(U)+f(\nabla w\right)]\right)&\textrm{in }\Omega\times (0,\infty)\\
    w_{tt}+\Delta^2w-\Delta w_{tt}=\textrm{div}\left(\mathcal{C}\left[e(U)+f(\nabla w\right)]\nabla w\right)&\textrm{in }\Omega\times (0,\infty),
\end{cases}\]
where $\mathcal{C}$ is a symmetrical $2\times 2$ matrix defined as
\[C[e]=\dfrac{E}{d(1-\mu^2)}[\mu(\mathrm{Tr} e)I+(1-\mu)e],\]
describes the deflections of a $2D$ plate occuppying the domain $\Omega,$ where $U=(u,v)$ represents the in-plane displacement of the plate, and $w$ represents the vertical displacement of the plate. This model has been studied by several authors; in particular by Ciarlet \cite{ciar}, Sánches \cite{Luyo}, Lasiecka \cite{lasi2}, Lasiecka-Benabdallah \cite{lasi3}, Koch-Lasiecka \cite{lasi4}, and Puel-Tucsnak \cite{puel}. In 2000, Perla and Zuazua \cite{perla2} showed that the much simpler \emph{Timoshenko's model} for the vertical displacement $w,$
\[w_{tt}+\Delta^2 w-\Delta w_{tt}-c\left(\int_{\Omega}|\nabla w|^2 \mathrm{d}A\right)\Delta w=0,\,\,\,\textrm{in }\Omega\times (0,\infty),\]
may be obtained as the weak$^{\ast}$ limit of the von Karmán system when some parameters tend to zero. The term $-\Delta w_{tt}$ appears to absorb the rotational inertia effects at the point $x$ of the plate. The term $(-\Delta)^{\theta} u_t$ represents a fractional dissipation of energy, which becomes clear in \eqref{eq_energia}. The Timoshenko model with rotational inertia and fractional dissipation then is written as
\begin{equation}
    u_{tt}+\Delta^2 u -\Delta u_{tt}+(-\Delta)^{\theta}u_t=0.
\end{equation}
In the model studied in this work, the choice of the memory nonlinearity term is due to theoretical reasons explained before, and the term $-\Delta u$ is added to avoid singularity problems at the origin. The existence of solutions for the problem without this term is more difficult and it is still an open problem. 

We are interested into find nonexistence results for global in-time solutions for the subcritical case. In \cite{longen}, global existence is proved for dimension $n\le 4.$ In particular, for $n\le 3,$ in the case $\theta\in \left[ 0,\frac12\right),$ $\gamma\in (0,1)$ and $p>\overline{p},$ the critical exponent is $\overline{p}=\max\left\lbrace p_c,\gamma^{-1}\right\rbrace,$ where $p_c=p_c(n,\gamma,\theta)$ is given by

\begin{equation}\label{pcrit}
p_c(n,\gamma,\theta):=1+\frac{2\left(1+(1-\gamma)(1-\theta)\right)}{\left(n-2+2\gamma(1-\theta)\right)_+},
\end{equation}

As mentioned before, this kind of competition between two different values, the first being a Fujita-type and the second being $\gamma^{-1},$ related to the nonlinear memory term, had already been noticed in precedent papers, for example \cite{DA14NODEA}, \cite{cazenave}. However, for dimension $n\ge 3$ a new factor comes into play. The possibility of $p<2,$ a scenario in which $L^p-L^q$ estimates are difficult to obtain, gives a worse decay rate. Even if for $n=3$ this loss in decay rate is not significant enough to produce a change in the critical exponent, for $n\ge 4$ the effect of this loss becomes apparent. In fact, if $\gamma$ is large enough in comparison to $\theta,$  a new exponent joins the competition,
\begin{equation}
    \label{p_c_til}
    \tilde{p}_c:=\dfrac{6n+4-4(n+1)\theta}{2s+n(3-2\theta)-4(1-\gamma)(1-\theta),}
\end{equation}
where $s$ is the regularity of initial data, $(u_0,u_1)\in H^s(\rn)\cap L^1(\rn)\times H^{s-1}(\rn)\cap L^1(\rn).$ Therefore, for $n\ge 4,$ the critical exponent becomes
\[\overline{p}=\max\lbrace\gamma^{-1},p_c,\tilde{p}_c\rbrace.\]

\begin{figure}[!h]
\begin{center}
\includegraphics[width=.48\linewidth]{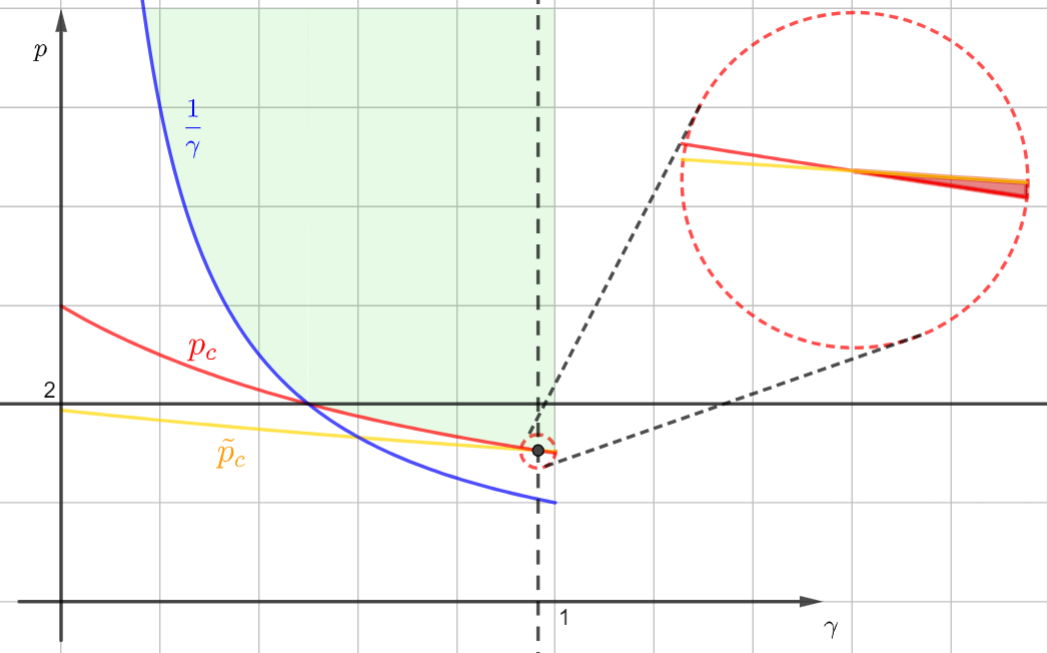}  
\caption{Critical Exponent for $n\ge4$}
\label{fig:pcrit}
\end{center}
\end{figure}

The reader should check Figure \ref{fig:pcrit} to see the small region where $\tilde{p}_c>p_c.$ The figure illustrates the case $n=4$ and $\theta=\frac{1}{100}.$

The goal of this work is to answer the following question: Is this new exponent optimal or is it just fruit of a suboptimal result for the decay rate when $p<2?$ If one can prove nonexistence for $p<\Tilde{p_c},$ then the problem is solved, with exception, possibly, to the critical case. Otherwise, if one does not find $p<\Tilde{p_c}$ in nonexistence results, then the region colored in red in Figure \ref{fig:pcrit} remains a mystery. To answer this question, we study the nonexistence of solutions in the subcritical case.

Our main result in this paper reads as follows.
\begin{teorema} \label{teorema}
Let $p_c$ be as in Definition \eqref{pcrit}, define
\[\overline{p}=\begin{cases}
\gamma^{-1},&,\mathrm{if}\,\,\gamma\in\left(0,\frac{n-2}{n}\right]\\
p_c,&\mathrm{if}\,\,\gamma\in\left(\frac{n-2}{n},1\right)
\end{cases}\]
and fix $q=n+2\theta.$ Assume that $u_0,u_1\in L^1(\langle x \rangle^q dx).$ Moreover, assume the sign condition
\[\IRn{u_1}>0.\]
If there exists a global in-time (nontrivial) weak solution
\[u\in L^p\left([0,\infty),L^p(\rn,\langle x\rangle^{-q})\right)\]
to problem (\ref{P}), then $p\geq\overline{p};$

Furthermore, in the case $\gamma\in\left(\frac{n-2}{n},1\right)$ the inequality is strict, that is, $p>\overline{p}=p_c.$
\end{teorema}

\section{Test Functions and Weak Solutions}

In this work, we find global in-time solutions to the problem (\ref{P}) for the supercritical case, for a specific \emph{critical exponent} candidate. The sharpness of such candidate is achieved by showing nonexistence of global solutions in the subcritical case.

This nonexistence counterpart is usually derived by using a classic test function method. However, since the fractional Laplace operator has a nonlocal behavior, such method is not applicable (at least directly), since it relies on the compactness of the test functions'  support under the action of these operators.

Having this in mind, we use a modified test function to obtain our results. To deal with the nonlocality of the fractional Laplace operators, we will replace the usual compactly supported test functions by some suitable test functions with polynomial decay. To this end, we introduce a class of functions and a definition of a weak solution to problem (\ref{P}) which is adequate to our purposes.
\begin{definicao}
Let $\theta\in(0,1)$ be a fixed number, and fix $q=n+2\theta.$ We define the space $\mathcal{C}_q^\infty(\rn)$ as the subspace of infinitely differentiable functions $\varphi$ such that $\left\langle x\right\rangle^q\varphi$ is bounded, and for any $\sigma>0$ with $\sigma$ integer or $\sigma-\lfloor\sigma\rfloor\in[\theta,1),$ the function $\left\langle x\right\rangle^q(-\Delta)^\sigma\varphi$ is also bounded.
\end{definicao}

The choice of this specific space $C^\infty_q(\rn),$ with $q=n+2\theta$ is justified by the following lemma and its corollary. With these, we may show that a function in $C^\infty_q(\rn)$ remains in the same space after the action of the Fractional Laplace operator $(-\Delta)^{\theta}.$ First, we recall the definition of the fractional Laplace operator and an alternative definition for non-integer powers, that is equivalent with the usual definition when the domain is the whole space $\rn.$

\begin{definicao}\label{def:fraclap}
For any $\sigma>0,$ we may define the fractional Laplace operator $(-\Delta)^\sigma: H^{2\sigma}\to L^2$ as
\begin{equation}
    (-\Delta)^\sigma f=\mathcal{F}^{-1}\left(|\xi|^{2\sigma}\hat{f}\right).
\end{equation}
If $\sigma$ is not integer, then the operator $(-\Delta)^\sigma$ admits an integral representation. For $y\in\rn,$ let $\tau_y$ be the translation operator, i.e.,
\begin{equation}
    \tau_yf(x)=f(x+y).
\end{equation}
Then the identity
\begin{equation}
    (-\Delta)^{\sigma}f(x)=(-1)^{\lfloor\sigma\rfloor+1}C_\sigma\int_{\rn}\frac{(\tau_{\frac{y}{2}}-\tau_{-\frac{y}{2}})^{2\lfloor\sigma\rfloor+2}f(x)}{|y|^{n+2\sigma}}~\d y
\end{equation}
holds for any $f\in H^{2\sigma}.$ Also, the constant $C_\sigma$ is given by
\[C_\sigma = 2^{-2\lfloor\sigma\rfloor-2+2\sigma} \left(\int_{\rn}\frac{\sin(y_1)^{2\lfloor\sigma\rfloor+2}}{|y|^{n+2\sigma}}~\d y\right)^{-1}>0.\]
\end{definicao}
\begin{observacao}
The identity given in Definition \ref{def:fraclap} is shown in \cite{nezza}.
\end{observacao}

The fractional Laplace operator can be conveniently extended to more general spaces. In particular, it may be extended by duality to the tempered distribution space $\mathcal{S}'.$ The following Lemma and its Corrolary ensure that a bounded function with bounded derivatives can be controlled under the action of the fractional Laplace operator.

\begin{lema}\label{lema:japCq}
Assume $f\in C^2$ bounded, with bounded derivatives. If there exists a constant $C_0$ such that the estimates
\[|f(y)|\leq C_0 |f(x)|,\,\,\,\,\sup_{|\alpha|=2}|\partial^\alpha f(y)|\leq C_0 \sup_{|\alpha|=2}|\partial^\alpha f(x)|\]
hold when $|x|\leq |y|,$ then for $|x|>1,$ the following pointwise estimate holds:
\begin{align}
    |(-\Delta)^{\sigma}f(x)|&\leq C |x|^{-n-2\sigma}\int_{|y|<3|x|}|f(y)|\d y+ C|f(x)||x|^{-2\sigma}\notag\\&+C|x|^{2-2\sigma}\sum_{|\alpha|=2}\frac{|\alpha|}{\alpha!}\big|\partial^\alpha f\left(\tfrac{x}{2}\right)\big|,
\end{align}
for any $\sigma\in(0,1).$
\end{lema}
\begin{prova} See \cite{fujiwara}.
\end{prova}

As a consequence of Lemma \ref{lema:japCq}, one can derive the following corollary, which bounds the action of the fractional Laplace operator on $\langle x \rangle^{-q}$ pointwisely.

\begin{corolario}\label{cor_Cinftyq}
Let $f(x)=\langle x \rangle^{-q}$, for $q>n,$ and let $\sigma>0.$ We set $s=\sigma-\lfloor\sigma\rfloor$. Then,
\[\big|(-\Delta)^\sigma f(x) \big|\leq C \langle x \rangle^{-q_\sigma},\qquad \forall~x\in\rn,\]
where $q_\sigma=q+2\sigma$ if $\sigma$ is an integer, or $q=n+2s$ otherwise.
\end{corolario}

\begin{prova}
We start with the following identity
\begin{align*}
-\Delta f(x)&=q\sum_{j=1}^n \partial_{x_j}\left(x_j\langle x \rangle^{-q-2}\right)\\
%&=q\sum_{j=1}^n\Big(\langle x \rangle^{-q-2}-(q+2)x_j^2\langle x \rangle^{-q-4} \Big)\\
%&=qn\langle x \rangle^{-q-2}-q(q+2)|x|^2\langle x \rangle^{-q-4}\\
%&=qn\langle x \rangle^{-q-2}-q(q+2)\langle x \rangle^{-q-2}+q(q+w)\langle x \rangle^{-q-4}\\
&=q(q+2-n)\langle x \rangle^{-q-2}+q(q+2)\langle x\rangle^{-q-4}\\
&:=c_0\langle x \rangle^{-q-2}+c_1\langle x\rangle^{-q-2-2},
\end{align*}
with $c_0,c_1$ depending only on $n,\sigma,q.$ After $\lfloor\sigma\rfloor$ iterations of this identity, we get
\begin{equation}\label{eq:lap_japq}
(-\Delta)^{\lfloor\sigma\rfloor}f(x)=\sum_{k=0}^{\lfloor\sigma\rfloor}c_k\langle x \rangle^{-q-2\lfloor\sigma\rfloor-2k},
\end{equation}
for some $c_k=c_k(n,\sigma,q)\in \R.$

Assuming that $\sigma$ is integer, then $\lfloor\sigma\rfloor=\sigma$ and $s=0.$ Therefore, \eqref{eq:lap_japq} implies
\[\Big|(-\Delta)^{\sigma}f(x) \Big|\leq C\langle x \rangle^{-q-2\sigma},\]
which is precisely what we needed to prove in this case. Next, assume $\sigma$ non-integer, that is, $s\in(0,1).$ We have
\begin{align}\label{cor:id}
\Big|(-\Delta)^{\sigma}f(x) \Big|&=\Big|(-\Delta)^{s}\left((-\Delta)^{\lfloor\sigma\rfloor}f(x)\right) \Big|\notag\\
&=\Bigg|(-\Delta)^{s}\sum_{k=0}^{\lfloor\sigma\rfloor}c_k\langle x \rangle^{-q-2\lfloor\sigma\rfloor-2k} \Bigg|\notag\\
&\leq \sum_{k=0}^{\lfloor\sigma\rfloor}|c_k |~ \Big| (-\Delta)^{s}\langle x \rangle^{-q-2\lfloor\sigma\rfloor-2k}\Big|.
\end{align}

Then, applying Lemma \ref{lema:japCq} for each $g(x)=\langle x \rangle^{-q-2\lfloor\sigma\rfloor-2k}$ with $k=0,... \lfloor\sigma\rfloor$ (plugging $s$ in the place of $\sigma$), we obtain the estimates
\begin{align*}
\Big| (-\Delta)^{s}\langle x \rangle^{-q-2\lfloor\sigma\rfloor-2k}\Big|&\leq C|x|^{-n-2s}\int_{|y|<3|x|}|g(y)|\d y+C|g(y)||x|^{-2s}\\&+C|x|^{-2s+2}\sum_{|\alpha|=2}\big| \partial^\alpha g\left(\tfrac x2\right)\big|,
\end{align*}
for every $|x|>1.$ Now, since $g$ and $\partial^{\alpha}g$  are bounded, we must show only that the integral on the right-hand side is bounded. But this is a consequence of $\langle x \rangle^{-a}\in L^1(\rn)$ for $a>n.$ Indeed,
\begin{align*}
\int_{|y|<3|x|}|g(y)|~\d y&\leq \int_{\rn} |g(y)|~\d y \\&\leq \int_{|y|<1}|g(y)|~\d y+\int_{|y|>1}|g(y)|~\d y\\
&\leq \int_{|y|<1}\d y + \int_{|y|>1}|y|^{-a}\d y<+\infty,~\text{for}~a>n.
\end{align*}
Since for $|x|>1$ we have $|x|\sim \langle x \rangle,$ we get the estimates
\begin{equation}\label{cor:sub}
\Big| (-\Delta)^{s}\langle x \rangle^{-q-2\lfloor\sigma\rfloor-2k}\Big|\leq C\langle x \rangle^{-n-2s},
\end{equation}
for $k=0,1,...\lfloor\sigma\rfloor.$
Therefore, using \eqref{cor:sub} in \eqref{cor:id}, we obtain
\begin{equation}
\Big|(-\Delta)^{\sigma}f(x) \Big|\leq C \langle x \rangle^{-n-2s},
\end{equation}
which concludes our proof.
\end{prova}

\begin{observacao}
\label{remarkphi}
We observe that, choosing $q=n+2\theta,$ with $\theta\in[0,1)$ Corollary \ref{cor_Cinftyq} says that $f(x)=\langle x \rangle^{-q}\in C_q^\infty(\rn),$ and that $(-\Delta)^{\sigma}\langle x \rangle^{-q}\in C_q^\infty(\rn)$ for every $\sigma>0.$
\end{observacao}

Having defined the appropriate function space to our problem, we may introduce a definition of a weak solution to problem \eqref{P}, locally and globally with respect to the time variable.

\begin{definicao}\label{defsolfraca}
Fix $q=n+2\theta,$ and fix $T\in[0,\infty).$ We say that $u\in L^p_{loc}\big([0,T),L^p\left(\rn,\left\langle x\right\rangle^q\d x\right)\big)$ is a weak solution to \eqref{P} if for any function $\psi\in \mathcal{C}^2_c([0,T])$ satisfying $\psi(0)= 1$, $\psi(T)=\psi_t(T)=0,$ and for any $\varphi\in\mathcal{C}_q^\infty(\rn),$ it holds
\small\begin{align} \label{weaksolution}
   \Iabt{0}{T}{\psi(t)\IRn{&F(t,u)\varphi(x)}}\\&= \Iabt{0}{T}{\psi(t)\IRn{u(t,x)(-\Delta\varphi)(x)}}+\Iabt{0}{T}{\psi(t)\IRn{u(t,x)(\Delta^2\varphi)(x)}}\notag\\
   &-\Iabt{0}{T}{\psi_t(t)\IRn{u(t,x)((-\Delta)^\theta\varphi)(x)}}+\Iabt{0}{T}{\psi_{tt}(t)\IRn{u(t,x)\varphi(x)}}\notag\\
   &+\Iabt{0}{T}{\psi_{tt}(t)\IRn{u(t,x)(-\Delta\varphi)(x)}}-\IRn{u_0(x)((-\Delta)^\theta \varphi)(x)}\notag\\
   &-\IRn{u_1(x) \varphi(x)}-\IRn{u_1(x)(-\Delta \varphi)(x)}\notag\\
    &+\psi_t(0)\left(\IRn{u_0(x)\varphi(x)}+\IRn{u_0(x)(-\Delta\varphi)(x)}\right).
\end{align}

\normalsize We say that the weak solution is locally-in-time defined if $T<\infty$ and is globally-in-time defined if $T=\infty.$

Equivalently, a function $u\in L^p_{loc}\big([0,T),L^p\left(\rn,\left\langle x\right\rangle^q\d x\right)\big)$ is a global in-time solution if, and only if, $u\big|_{[0,T)\times \rn}$ is a local in-time weak solution, for any $T>0.$
\end{definicao}

\begin{observacao}
We remark that it holds
\[L^p_{loc}\big([0,T),L^p\left(\rn,\left\langle x\right\rangle^q\d x\right)\big)\subset L^p_{loc}\big([0,T),L^p\left(\rn\right)\big), \]
so the above defined weak solutions space is properly contained in a more conventional solution space.
\end{observacao}

We may show with ease that classical solutions to problem (\ref{P}) are also weak solutions.

\begin{prop}
Assume that $u_0,u_1 \in \mathcal{S}(\rn).$ Also, assume that $u\in \mathcal{C}^2\left([0,T),\mathcal{S}\right)$ is a ``classical'' solution to problem (\ref{P}). Then, $u$ is also a weak solution to (\ref{P}), according to Definition \ref{defsolfraca}.
\end{prop}
\begin{prova}
Multiplying the equation in (\ref{P}) by $\psi(t)\varphi(x)$ and integrating in $[0,T]\times \rn,$ we get
\small\begin{align}
    \Iabt{0}{T}{\psi(t)\IRn{&F(t,u)\varphi(x)}}\notag \\&= \Iabt{0}{T}{\psi(t)\IRn{(-\Delta  u)(t,x)\varphi(x)}}\notag +\Iabt{0}{T}{\psi(t)\IRn{(\Delta^2u)(t,x)\varphi(x)}}\\
    & +\Iabt{0}{T}{\psi(t)\IRn{((-\Delta)^\theta u_t)(t,x)\varphi(x)}} +\Iabt{0}{T}{\psi(t)\IRn{u_{tt}(t,x)\varphi(x)}}\notag\\
    & +\Iabt{0}{T}{\psi(t)\IRn{(-\Delta u_{tt})(t,x)\varphi(x)}}.
\end{align}

\normalsize Integrating by parts in space, due to $u,u_t, u_{tt}\in \mathcal{S}$ and $\varphi\in \mathcal{C}^\infty_c(\rn),$ 
\begin{align}
    \Iabt{0}{T}{\psi(t)\IRn{&F(t,u)\varphi(x)}}\notag \\&= \Iabt{0}{T}{\psi(t)\IRn{u(t,x)(-\Delta\varphi)(x)}}+\Iabt{0}{T}{\psi(t)\IRn{u(t,x)(\Delta^2\varphi)(x)}}\notag\\
    & +\Iabt{0}{T}{\psi(t)\IRn{ u_t(t,x)((-\Delta)^\theta\varphi)(x)}}+\Iabt{0}{T}{\psi(t)\IRn{u_{tt}(t,x)\varphi(x)}}\notag\\
    & +\Iabt{0}{T}{\psi(t)\IRn{u_{tt}(t,x)(-\Delta\varphi)(x)}}.
\end{align}

Now, integrating by parts in time as many times as needed to get rid of all time derivatives of $u$, recalling that $\psi(0)= 1$, $\psi(T)=\psi_t(T)=0,$ we obtain
\small\begin{align}
    \Iabt{0}{T}{&\psi(t)\IRn{F(t,u)\varphi(x)}}\notag\\&= \Iabt{0}{T}{\psi(t)\IRn{u(t,x)(-\Delta\varphi)(x)}}+\Iabt{0}{T}{\psi(t)\IRn{u(t,x)(\Delta^2\varphi)(x)}}\notag\\
    &-\IRn{u(0,x)((-\Delta)^\theta \varphi)(x)}-\Iabt{0}{T}{\psi_t(t)\IRn{ u(t,x)(-\Delta)^\theta\varphi(x)}}\notag\\
    &-\IRn{u_t(0,x)\varphi(x)}+\IRn{\psi_t(0)u(0,x)\varphi(x)}\notag\\&+\Iabt{0}{T}{\psi_{tt}(t)\IRn{u(t,x)\varphi(x)}}-\IRn{u_t(0,x)(-\Delta)\varphi(x)}\notag\\
    &+\IRn{\psi_t(0)u(0,x)(-\Delta \varphi)(x)}+\Iabt{0}{T}{\psi_{tt}(t)\IRn{u(t,x)(-\Delta)\varphi(x)}}.
\end{align}

\normalsize Lastly, applying the boundary conditions $u(0,x)=u_0(x)$ and 
$u_t(0,x)=u_1(x)$ we arrive exactly at (\ref{weaksolution}), thus concluding the proof.

\end{prova}

\vspace{20pt}
The idea of using polynomially decaying test functions instead of the usual compactly supported test functions is quite recent, from D'Abbicco and Fujiwara in 2021 \cite{fujiwara}. In the problem addressed in \cite{fujiwara}, the nonlinearity is a $p-$power time-derivative of $u,$ that is, $|\partial_t^\ell u|^p,$ with $\ell$ positive integer. Since in the present work the nonlinearity is a memory term, we must slightly change the argument, specifically on the time-related part. To this end, we introduce the fractional and differential operators as follows, as well as an adequate time-dependent function. For the fractional and differential operators, we may follow the definition given in \cite{samko}. We refer the reader to \cite{samko} for more properties related to these operators.

\begin{definicao}
Let $\alpha\in (0,1)$ and fix $~T>0.$ We define, for any function $f\in L^1(0,T),$

\end{definicao}
\label{def:J}
\begin{equation}
    J_{0|t}^\alpha f(t):=\frac{1}{\Gamma(\alpha)}\Iot{(t-s)^{-(1-\alpha)}f(s)},
\end{equation}
\begin{equation}
    J_{t|T}^\alpha f(t):=\frac{1}{\Gamma(\alpha)}\Iab{t}{T}{(s-t)^{-(1-\alpha)}f(s)},
\end{equation}
called the left-sided and right-sided Riemann-Liouville fractional integrals of the order $\alpha.$

\begin{definicao}
Let $\alpha\in (0,1)$ and fix $~T>0.$ We define, for any function $f\in AC([0,T]),$ the space of all absolutely continuous functions on $[0,T],$

\begin{equation}
    D_{0|t}^\alpha f(t) := \partial_t J_{0|t}^{1-\alpha}f(t),
\end{equation}
\begin{equation}
     D_{t|T}^\alpha f(t) := -\partial_t J_{t|T}^{1-\alpha}f(t),
\end{equation}
which are called the left-sided and right-sided Riemann-Liouville fractional derivatives.
\end{definicao}

\begin{definicao}
\label{defJD}
Let $\alpha\in(0,1)$ and fix $T>0.$ We define, for any function $f\in L^1(0,T),$
\begin{align*}
    J_{0|t}^\alpha f(t):=\frac{1}{\Gamma(\alpha)}\Iot{(t-s)^{-(1-\alpha)}f(s)},\,\,\,\,\,\,\,\,\,\,& D_{0|t}^\alpha := \partial_t J_{0|t}^{1-\alpha}\\
     J_{t|T}^\alpha f(t):=\frac{1}{\Gamma(\alpha)}\Iab{t}{T}{(s-t)^{-(1-\alpha)}f(s)},\,\,\,\,\,\,\,\,\,\,& D_{t|T}^\alpha := -\partial_t J_{t|T}^{1-\alpha},
\end{align*}

\end{definicao}

The following theorem asserts that the Riemann-Liouville Integral and Derivative are, in some sense, inverse of each other. For its proof and a more detailed explanation, the reader is addressed to \cite{samko}.

\begin{teorema}[\cite{samko}, p.44] \label{propJD1} Let $\alpha\in(0,1).$ Then the equality
\begin{equation}
D_{0|t}^\alpha J_{0|t}^\alpha f(t)= f(t)
\end{equation}
is valid for any $f\in L^1([0,T]).$, Also, the equality
\begin{equation}\label{JD}
J_{0|t}^\alpha D_{0|t}^\alpha f(t)= f(t)
\end{equation}
is satisfied for $f(t)\in J_{0|t}^\alpha(L^1([0,T])).$
\end{teorema}

With this essential property, we can also prove a generalization of the integration by parts formula.

\begin{prop}[Integration by Parts Formula] Assume that $\alpha\in(0,1),~\varphi\in L^p([0,T]),~\psi\in L^q([0,T]),$ and $\frac1p+\frac1q\leq 1+\alpha.$ Then, the formula for fractional integration by parts,
\begin{equation}\label{JDparts}
\Iabt{0}{T}{\varphi(t)(J_{0|t}^\alpha \psi)(t)}=\Iabt{0}{T}{\psi(t)(J_{t|T}^\alpha \varphi)(t)}
\end{equation}
is valid.
\end{prop}
\begin{prova}
The proof is based in interchanging the order of integration, by making use of Fubini's Theorem. To justify its application, we argue as follows. First, from the Hardy-Littlewood Theorem with limiting exponent (see \cite{samko},Th.3.5), the fractional integration operator $J_{0|t}^\alpha$ is bounded from $L^p([0,T])$ into $L^q([0,T]),$ with $q=\frac{p}{1-\alpha p}.$ Therefore, one can apply Hölder inequality in \eqref{JDparts} and see that both integrals are absolutely convergent, and thus we can use Fubini's Theorem. We have then
\begin{align*}
\Iabt{0}{T}{\varphi(t)(J_{0|t}^\alpha \psi)(t)}&=\Iabt{0}{T}{\varphi(t)\frac{1}{\Gamma(\alpha)}\int_0^t(t-s)^{\alpha-1}\psi(s)\d s}\\
&=\frac{1}{\Gamma(\alpha)}\int_0^T \int_0^T \varphi(t)\chi_{[0,t]}(s)\, (t-s)^{\alpha-1}\psi(s)~\d s \d t\\
&=\frac{1}{\Gamma(\alpha)}\int_0^T \psi(s)\int_0^T \varphi(t)\chi_{[s,T]}(t)\, (t-s)^{\alpha-1}~\d t \d s\\
&=\frac{1}{\Gamma(\alpha)}\int_0^T \psi(s)\int_s^T \varphi(t)(t-s)^{\alpha-1}~\d t \d s\\
&=\int_0^T \psi(s)(J_{t|T}^\alpha\varphi)(s)~\d s,
\end{align*}
where $\chi_A(x)$ represents the indicator function of $x\in A,$ that is, $\chi_A(x)=1,$ if $x\in A$ and $\chi_A(x)=0, $ if $x\notin A.$
 \end{prova}

As a consequence, we can prove also a version of the integration by parts formula with Riemann-Liouville Integrals in the argument of integration over $[0,T].$

\begin{prop}\label{propJD2}
Assume that $\alpha\in(0,1),~f\in J_{t|T}^\alpha (L^p([0,T])), ~g\in J_{t|T}^\alpha (L^q([0,T]))$ with $\frac1p+\frac1q\leq 1+\alpha.$ Then, the formula
\begin{equation}\Iabt{0}{T}{(D_{0|t}^\alpha f)(t)g(t)}=\Iabt{0}{T}{f(t)(D_{t|T}^\alpha g)(t)}\end{equation}
holds.
\end{prop}

\begin{prova}
Let $\varphi(t)=D_{0|t}^\alpha f(t)\in L^p([0,T])$ and $\psi(t)=D_{t|T}^\alpha g(t)\in L^q([0,T]).$ Then, using \eqref{JD} and integration by parts, we get
\begin{align*}
\Iabt{0}{T}{(D_{0|t}^\alpha f)(t)g(t)}&=\Iabt{0}{T}{\varphi(t) g(t)}=\Iabt{0}{T}{\varphi(t)\left(J_{t|T}^\alpha D_{t|T}^\alpha g\right)(t)}\\
&=\Iabt{0}{T}{\varphi(t)(J_{t|T}^\alpha \psi)(t)}
=\Iabt{0}{T}{(J_{0|t}^\alpha \varphi)(t)\psi(t)}\\
&=\Iabt{0}{T}{\left(J_{0|t}^\alpha D_{0|t}^\alpha f\right)(t)\psi(t)}
=\Iabt{0}{T}{f(t)(D_{t|T}^\alpha g(t))}.
\end{align*}

\end{prova}

%\begin{prop}\label{propJD}
%The operators defined in Definition \ref{defJD} satisfy the following properties:
%\begin{enumerate}[label=\roman*)]
 %   \item \begin{equation}\Iabt{0}{T}{(D_{0|t}^\alpha f)(t)g(t)}=\Iabt{0}{T}{f(t)(D_{t|T}^\alpha g)(t)}\end{equation}
%    \item \begin{equation}D_{0|t}^\alpha J_{0|t}^\alpha f(t)= f(t). \end{equation}
%\end{enumerate}
%\end{prop}

\begin{definicao}\label{defomega}
We define the following auxiliary function, for a fixed $T>0:$
\begin{equation}
    \omega_T(t)=\begin{cases}
    1-\frac{t}{T},& t\in[0,T]\\
    0,& t>T.
    \end{cases}
\end{equation}
\end{definicao}

\begin{observacao}
We remark that $\supp{\omega_T}=[0,T]$ and that $(\omega_T(t))^\beta \in \mathcal{C}_c^k(\R),$ for any $\beta>k.$ It also has the following interaction with the fractional derivative, which is of our interest:
\end{observacao}

\begin{lema}
\label{lemapsiT}
For any $\alpha\in(0,1),$ it follows that

\begin{equation}\label{eq:Domega}
    D_{t|T}^\alpha \omega_T(t)^\beta = C_{\alpha,\beta} T^{-\alpha} \omega_T(t)^{\beta-\alpha}, 
\end{equation}

where
\[C_{\alpha,\beta}= \frac{\Gamma(\beta+1)}{(\beta+2-\alpha)\Gamma(\beta-\alpha)}.\]
\end{lema}
\begin{prova}
We have initially that, making the  change of variables $\tau= \frac{s-t}{T-t},$
\begin{align*}
  \Gamma(1-\alpha)J_{t|T}^{1-\alpha} \omega_T(t)^\beta &=\Iab{t}{T}{(s-t)^{-\alpha}\left(1-\tfrac{s}{T}\right)^\beta}\\&=\Iabtau{0}{1}{\tau^{-\alpha}(T-t)^{-\alpha}T^{-\beta}(T-t)^\beta (1-\tau)^\beta (T-t)}\\
  &= (T-t)^{-\alpha+\beta+1}T^{-\beta} \Iabtau{0}{1}{\tau^{-\alpha}(1-\tau)^\beta}\\
  &= \omega_T(t)^{-\alpha+\beta+1} T^{1-\alpha}\,\,\,\frac{\Gamma(1-\alpha)\Gamma(1+\beta)}{\Gamma(\beta+2-\alpha)}.
\end{align*}

Hence,
\begin{equation}
\label{eqJomega}
    J_{t|T}^{1-\alpha} \omega_T(t)^\beta= \frac{\Gamma(1+\beta)}{\Gamma(\beta+2-\alpha)}\,T^{1-\alpha}\,\,\omega_T(t)^{-\alpha+\beta+1}.
\end{equation}

We now differentiate with respect to $t,$ obtaining
\begin{equation}\label{eqDtomega}
  \partial_t \omega_T(t)^{-\alpha+\beta+1}=(-\alpha+\beta+1)\left(1-\tfrac{t}{T}\right)^{-\alpha+\beta}(-T^{-1}).  
\end{equation}

Using (\ref{eqDtomega}) in (\ref{eqJomega}), we get
\begin{align*}
    D_{t|T}^\alpha \omega_T(t)^\beta &= \frac{\Gamma(1+\beta)}{\Gamma(\beta+2-\alpha)}\,T^{1-\alpha}(-\alpha+\beta+1)\left(1-\tfrac{t}{T}\right)^{-\alpha+\beta}(-T^{-1})\\
    &=\frac{\Gamma(\beta+1)}{(\beta+2-\alpha)\Gamma(\beta-\alpha)}\,\,T^{-\alpha}\,\omega_T(t)^{\beta-\alpha}.
\end{align*}
\end{prova}

Lastly, differentiating twice the identity from Lemma \ref{lemapsiT} above, one obtains the following estimate, which will be useful in proving the main Theorem.

\begin{lema} For any $\alpha\in(0,1),$ and $j=0,1,2,$ the following holds:
    \begin{align}|\partial^j_t\psi_T(t)|&\leq C_{\alpha,\beta}T^{-\alpha-j}\omega_T(t)^{\beta-\alpha-j},\label{lemaDpsiT}
\end{align}
\label{lema:delTpsi}
\end{lema}

\vspace{10pt}

We conclude this section remarking that the function $D_{t|T}^\alpha \omega_T(t)^\beta$ satisfies all the requirements to be chosen as a time-dependent test function according to Definition \ref{defsolfraca}, for an appropriately chosen $\beta.$ Indeed, $\supp{D_{t|T}^\alpha \omega_T(t)^\beta}=[0,T],$ and it is in $\mathcal{C}^2,$ for any $\beta$ such that $\beta-\alpha>2,$ or equivalently, $\beta>2+\alpha.$ Also, from the last lemma, we have $D_{t|T}^\alpha \omega_T^\beta(T)=0,$ and the same holds for its derivative, that is, $\partial_t D_{t|T}^\alpha \omega_T^\beta(T)=0,$

\section{Nonexistence Results}

Since the Fractional Laplacian is a nonlocal operator, it is not easy, in general, to prove the optimality of $\overline{p}$ by applying the test function method. This section is dedicated to prove the nonexistence theorem for the subcritical case.

We start with a Lemma for estimating integrals of a specific kind. This integral appears recurrently in calculations of estimates in the main Theorem. 

\begin{lema}
\label{lema:Iu}
Let $\alpha\in(0,1), \beta>(\alpha+j)p,$ $\sigma,T>0, p,R>1,$ and $j=0,1,2.$ Let $\omega_T$ be as in Definition \ref{defomega}.

Assume that $\psi\in C^2(\R)$ and $\varphi\in C^\infty(\rn)$ satisfy the estimates
\[\big|\partial_t^j \psi\big|\le CT^{-\alpha-j}\omega_T^{\beta-\alpha-j}, \qquad \big|(-\Delta)^\sigma\varphi \big|\le CR^{-2\sigma}|\varphi|\].

Then, for any $\varepsilon>0,$
\begin{align*}
    \Bigg| \Iabt{0}{T}{\partial^j_t\psi&\IBR{u((-\Delta)^\sigma\varphi)}} \Bigg|\leq \varepsilon \Iabt{0}{T}{\IBR{\omega_T(t)^\beta |u|^p|\varphi(x)| }} + C_{\varepsilon,n} T^{1-(\alpha+j) p'}R^{n-2\sigma p'},
\end{align*}
for some constant $C_{n,\varepsilon}$ depending only on $n$ and $\varepsilon.$
\end{lema}

\begin{prova}
    Let $\varepsilon>0.$ Applying Hölder's Inequality,
\small\begin{align}
    \Bigg| \Iabt{0}{T}{\partial_t^j\psi(t)\IBR{u&(t,x)((-\Delta)^{\sigma}\varphi)(x)}} \Bigg|\notag \\& \leq  \varepsilon \Iabt{0}{T}{\IBR{\omega_T(t)^\beta |u|^p|\varphi(x)| }}\notag\\& +C_\varepsilon T^{-(\alpha+j)p'}R^{-2\sigma p'}\Iabt{0}{T}{\IBR{\omega_T^{\beta-(\alpha+j)p'}  |\varphi(x)|  }}\notag 
\end{align}
\normalsize Now since both $\omega_T$ and $\phi$ are bounded, and $\beta>(\alpha+j)p',$

\small\begin{align}
\Bigg| \Iabt{0}{T}{\partial_t^j\psi(t)\IBR{u&(t,x)((-\Delta)^{\sigma}\varphi)(x)}} \Bigg|\notag \\& \leq \varepsilon \Iabt{0}{T}{\IBR{\omega_T(t)^\beta |u|^p|\varphi(x)| }} + C_nC_\varepsilon T^{1-(\alpha+j)p'} R^{n-2\sigma p'}.
\end{align}
\normalsize Here, $C_n$ depends on the volume of the $n-$dimensional unit ball and the constant bounding $\phi.$
\end{prova}

We are finally ready to prove the main Theorem of this paper.

\vspace{15pt}
\textbf{Proof of Theorem \ref{teorema}:}

\begin{prova}[]
    We assume $p< \overline{p},$ and that $u$ is a nontrivial global in-time weak solution to (\ref{P}). We put $\alpha= 1-\gamma,$ fix $\beta>(\alpha+2)p'$ and fix suitable test functions depending on a parameter $R\gg 1$ as follows. 

%Let $\phi\in C_c^\infty$ be a smooth function satisfying
%\begin{itemize}
    %\item $\phi\big|_{B{\frac{1}{2}}}\equiv 1$;
    %\item $\supp\phi\subset\subset B_1$;
    %\item $\phi$ is radially decreasing.
%\end{itemize}
Consider the function $\omega_T$ as in Definition \ref{defomega}, and let $\Phi\in \mathcal{C}^\infty$ be defined as $\Phi(x):=\langle x\rangle^{-q}.$

We recall Remark \ref{remarkphi}, which states $\Phi(x)\in C_q^\infty(\rn)$ and $(-\Delta)^\theta\Phi(x)\in C_q^\infty(\rn).$
\begin{figure}[h]
    \centering
    \includegraphics[width=0.8\textwidth]{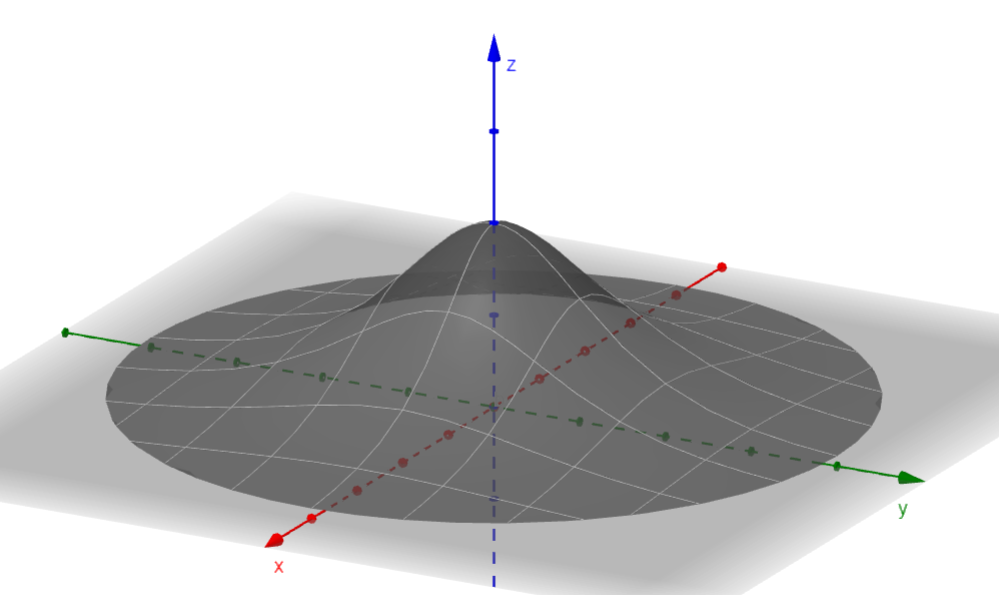}
    
    \caption{$\Phi\in C_q^\infty(\rn)$}
    \label{fig:enter-label}
\end{figure}

 For any $R\gg 1,$ and for a fixed $T> 0,$ we define
\[\varphi_R(x)=\Phi(R^{-1}x), \,\,\,\,\,\,\,\psi_T(t)=D_{t|T}^\alpha\left(\omega_T(t)^\beta\right).\]

First of all, we observe that, for any $\sigma\geq 0,$ 
\[\left((-\Delta)^\sigma \varphi_R\right)(x)\leq R^{-2\sigma}(-\Delta)^{\sigma}\varphi(x).\]

Since $u$ is a global weak solution to \eqref{P}, it satisfies
\small\begin{align}
   \Iabt{0}{\infty}{&\psi_T(t)\IRn{F(t,u)\varphi_R(x)}}\notag\\&= \Iabt{0}{\infty}{\psi_T(t)\IRn{u(t,x)(-\Delta\varphi_R)(x)}}\notag+\Iabt{0}{\infty}{\psi_T(t)\IRn{u(t,x)(\Delta^2\varphi_R)(x)}}\notag\\
   &-\Iabt{0}{\infty}{\partial_t\psi_T(t)\IRn{u(t,x)((-\Delta)^\theta\varphi_R)(x)}}+\Iabt{0}{\infty}{\partial_{tt}\psi_T(t)\IRn{u(t,x)\varphi_R(x)}}\notag\\ &+\Iabt{0}{\infty}{\partial_{tt}\psi_T(t)\IRn{u(t,x)(-\Delta\varphi_R)(x)}}-\IRn{u_0(x)(-\Delta)^\theta \varphi(x)}\notag\\ &-\IRn{u_1(x) \varphi(x)}-\IRn{u_1(x)(-\Delta) \varphi(x)}.
\end{align}

\normalsize Firstly, we obtain the following identity: Since $u_0,u_1\in L^1(\langle x \rangle^q dx)$ and $\varphi_R, (-\Delta)^{\theta}\varphi_R, \Delta \varphi_R\in L^\infty (\langle x \rangle^{-q} dx),$ Lebesgue's Dominant Convergence Theorem ensures that
\begin{align}
    \lim_{R\to\infty}\Bigg(\IRn{u_0&((-\Delta)^\theta \varphi_R)} +\IRn{u_1 \varphi_R}+\IRn{u_1(-\Delta \varphi_R)}\Bigg)\notag\\&=\IRn{u_0\lim_{R\to\infty}((-\Delta)^\theta \varphi_R)} \notag\\ & +\IRn{u_1 \lim_{R\to\infty}\varphi_R}+\IRn{u_1\lim_{R\to\infty}(-\Delta \varphi_R)}\notag\\
    &=\IRn{u_0\lim_{R\to\infty}R^{-2\theta}(-\Delta)^\theta \varphi} +\IRn{u_1 \lim_{R\to\infty}\varphi_R}\notag\\&+\IRn{u_1\lim_{R\to\infty}R^{-2}(-\Delta) \varphi}\notag\\&=\IRn{u_1}.
\end{align}

Due to the sign assumption from hypothesis, the latter is positive. Hence, for $R$ large enough, 
\small\begin{align} \label{eq:est5I}
   \Iabt{0}{T}{\psi_T(t)\IRn{F(t,u)\varphi_R(x)}}&< \Iabt{0}{T}{\psi_T(t)\IRn{u(t,x)(-\Delta\varphi_R)(x)}}\notag\\&+\Iabt{0}{T}{\psi_T(t)\IRn{u(t,x)(\Delta^2\varphi_R)(x)}}\notag\\
   &-\Iabt{0}{T}{\partial_t\psi_T(t)\IRn{u(t,x)((-\Delta)^\theta\varphi_R)(x)}}\notag\\&+\Iabt{0}{T}{\partial_{tt}\psi_T(t)\IRn{u(t,x)\varphi_R(x)}}\notag\\ &+\Iabt{0}{T}{\partial_{tt}\psi_T(t)\IRn{u(t,x)(-\Delta\varphi_R)(x)}}.
\end{align}

\normalsize Now, from the properties from Propositions \ref{propJD1} and \ref{propJD2}, since  $F(t,u)=\Gamma(\alpha)J_{0|t}^\alpha |u|^p,$ the left-hand side can be seen as
\small\begin{align*}
    \Iabt{0}{T}{\psi_T(t)\IRn{F(t,u)\varphi_R(x)}}&=\Gamma(\alpha)\Iabt{0}{T}{\IRn{\omega_T(t)^\beta |u(t,x)|^p\varphi_R(x)}}\\&:=\Gamma(\alpha) I(u).
\end{align*} 

\normalsize On the other hand, we estimate the five terms in the right-hand side. To do so, we integrate in $[0,T]$ in time and in $B_R$ in space, and after that we control the speeds at which we let $T$ and $R$ go to infinity.

We recall here the following estimates:
\[|\partial_t^j\psi_T(t)|\leq C T^{-\alpha-j} \omega_T(t)^{\beta-\alpha-j}, \,j=0,1,2,\qquad |((-\Delta)^\sigma\varphi_R)(x)|\leq C R^{-2\sigma}|\varphi_R(x)|.\]

Hence, we can apply Lemma \ref{lema:Iu} to bound all the five integrals, obtaining:

\begin{itemize}
    \item  \begin{equation} \label{est:I1}\Bigg| \displaystyle\Iabt{0}{T}{\psi_T(t)\IBR{u(t,x)(-\Delta\varphi_R)(x)}} \Bigg|\leq \varepsilon I(u) + C_\varepsilon T^{1-\alpha p'}R^{n-2p'}
    \end{equation}
    \item \begin{equation} \label{est:I2}\Bigg| \displaystyle\Iabt{0}{T}{\psi_T(t)\IRn{u(t,x)(\Delta^2\varphi_R)(x)}} \Bigg|\leq \varepsilon I(u) + C_\varepsilon T^{1-\alpha p'}R^{n-4p'}\end{equation}
    \item \begin{equation} \label{est:I3}\Bigg| \displaystyle\Iabt{0}{T}{\partial_t\psi_T\IRn{u((-\Delta)^\theta\varphi_R)}} \Bigg|\leq \varepsilon I(u) + C_\varepsilon T^{1-(\alpha+1) p'}R^{n-2\theta p'}\end{equation}
    \item \begin{equation}\label{est:I4}\Bigg| \displaystyle\Iabt{0}{T}{\partial_{tt}\psi_T\IRn{u\varphi_R}} \Bigg|\leq \varepsilon I(u) + C_\varepsilon T^{1-(\alpha+2) p'}R^{n}\end{equation}
    \item \begin{equation}\label{est:I5}\Bigg| \displaystyle\Iabt{0}{T}{\partial_{tt}\psi_T\IRn{u(-\Delta\varphi_R)}} \Bigg|\leq \varepsilon I(u) + C_\varepsilon T^{1-(\alpha+2) p'}R^{n-2p'}\end{equation}
\end{itemize}

\normalsize Now, using \eqref{est:I1}, \eqref{est:I2}, \eqref{est:I3},\eqref{est:I4} and \eqref{est:I5} in \eqref{eq:est5I}, we find
\small\begin{align}
    (\Gamma(\alpha)-5\varepsilon)& \Iabt{0}{T}{\IBR{\omega_T(t)^\beta |u(t,x)|^p\varphi_R(x)}} \notag \\ &<C_\varepsilon\Big(T^{1-\alpha p'}R^{n-2p'}+T^{1-\alpha p'}R^{n-4p'}\notag\\& +T^{1-(\alpha+1) p'}R^{n-2\theta p'}+T^{1-(\alpha+2) p'}R^{n}+T^{1-(\alpha+2) p'}R^{n-2p'}\Big).
\end{align}

\normalsize Choosing $\varepsilon<\frac{1}{5},$ we get, for $R,T$ large enough,
\small\begin{align}\label{est:nonexistence}
     \Iabt{0}{T}{\IBR{\omega_T(t)^\beta& |u(t,x)|^p\varphi_R(x)}}\notag \\&< C_{\alpha,\varepsilon}\Big(T^{1-\alpha p'}R^{n-2p'}+T^{1-\alpha p'}R^{n-4p'}\notag\\&\hspace{20pt}+T^{1-(\alpha+1) p'}R^{n-2\theta p'}+T^{1-(\alpha+2) p'}R^{n}+T^{1-(\alpha+2) p'}R^{n-2p'}\Big).
\end{align}
\normalsize Now, we divide into three cases:
\begin{enumerate}
    \item $\gamma>\frac{n-2}{n}$ and $p<p_c;$
    \item $\gamma>\frac{n-2}{n}$ and $p=p_c,$ the critical case;
    \item $\gamma\le\frac{n-2}{n}$ and $p<\frac{1}{\gamma};$
\end{enumerate}

\vspace{10pt}
\begin{enumerate}
    \item 
First consider the case where $\gamma>\frac{n-2}{n}.$ We stress that this is always true for $n=1,2.$ For $\eta\geq 0,$ we set $T=R^\eta.$ Then, we have
\small\begin{align}\label{eq:TReta}
     \Iabt{0}{R^\eta}{\IBR{\omega_T(t)^\beta |u(t,x)|^p\varphi_R(x)}}\notag &< C_{\alpha,\varepsilon}\Big(R^{-(\alpha\eta+2) p'+ n+\eta}+ R^{-(\alpha\eta+4) p' n+\eta}\notag\\&+ R^{-((\alpha+1)\eta+2\theta)p'+n+\eta}+R^{-(\alpha+2)\eta p'+n+\eta}\notag \\ &+R^{-((\alpha+2)\eta+2)p'+n+\eta}\Big).
\end{align}

\normalsize Now, we define, for $\eta\geq0,$
\begin{align} \label{def:g}g(\eta)&=\min\big\lbrace \alpha\eta+2,\alpha\eta+4,(\alpha+1)\eta+2\theta,(\alpha+2)\eta,(\alpha+2)\eta+2 \big\rbrace\notag \\
&=\min\big\lbrace \alpha\eta+2,(\alpha+1)\eta+2\theta,(\alpha+2)\eta \big\rbrace\notag\\
&=\begin{cases}
(3-\gamma)\eta,&\mathrm{if\,\,\,}\eta\in[0,2\theta]\\
(2-\gamma)\eta+2\theta,&\mathrm{if\,\,\,}\eta\in(2\theta,2(1-\theta)],\\
(1-\gamma)\eta +2,&\mathrm{if\,\,\,} \eta\in(2(1-\theta),\infty).
\end{cases}
\end{align}
where in the last step we used that $\alpha=1-\gamma.$ This function $g(\eta)$ gives the fastest possible decay for which we can control all five terms in \eqref{eq:TReta}. Below, we see a graphical representation of it.
\begin{figure}[!h]
    \centering
    \label{fig:geta}\includegraphics[scale=0.55]{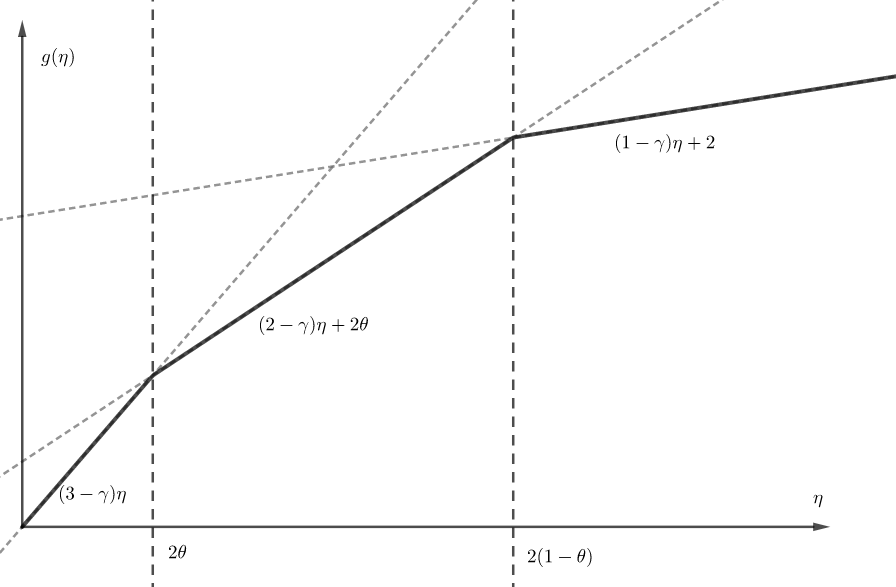}
    \caption{Function $g(\eta)$}
\end{figure}

Therefore, we have that
\begin{align}\label{estfinal}
     \Iabt{0}{R^\eta}{\IBR{\omega_T(t)^\beta |u(t,x)|^p\varphi_R(x)}}&< 5C_{\alpha,\varepsilon}R^{-g(\eta) p'+ n+\eta}.
\end{align}

To conclude the argument, it is sufficient that the exponent above is negative. This will happen if
\begin{align}
    g(\eta)p'-(n+\eta)>0,
\end{align}
or equivalently,
\begin{align}
    p'>\frac{n+\eta}{g(\eta)}\iff  p<\frac{n+\eta}{n+\eta-g(\eta)}=1+\frac{g(\eta)}{n+\eta-g(\eta)}:=h(\eta).
\end{align}

In particular, for $\eta=2(1-\theta),$ we have
\begin{align} h(2(1-\theta))&=1+\frac{2(1-\gamma)(1-\theta)+2}{n+2(1-\theta)-2(1-\gamma)(1-\theta)-2}\notag\\
&=1+\frac{2+2(1-\gamma)(1-\theta)}{n-2+2\gamma(1-\theta)}=p_c.
\end{align}

Moreover, we can see that $p_c$ is the greatest value for which we can prove nonexistence in this case. In other words, $h(\eta)$ reaches its maximum in $\gamma=2(1-\theta),$ assuming that $\gamma>\frac{n-2}{n}.$

Indeed, since
\[h'(\eta)=\frac{-g(\eta)+(n+\eta)g'(\eta)}{(n+\eta-g(\eta))^2},\]
its sign is determined by the sign of its numerator $-g(\eta)+(n+\eta)g'(\eta),$ which is given by
\[\mathrm{sgn}(-g(\eta)+(n+\eta)g'(\eta))=\begin{cases}
n(3-\gamma)>0,&\,\,\mathrm{if\,\,}\eta\in(0,2\theta)\\
(2-\gamma)-2\theta>0&\,\,\mathrm{if\,\,}\eta\in[2\theta,2(1-\theta))\\
-2+(1-\gamma)n,&\,\,\mathrm{if\,\,}\eta\in(2(1-\theta),\infty),
\end{cases}\]
and in the last case, we have
\[-2+(1-\gamma)n<0\iff \gamma>\frac{n-2}{n}.\]
Since $h(\eta)$ is continuous by parts, it has a local maximum point at $\eta=2(1-\theta)$ when $\gamma>\frac{n-2}{n}.$

Returning to \eqref{estfinal}, if $p<p_c,$ then applying Beppo-Levi's Monotone Convergence Theorem, since $\omega_T(t),\varphi_R(x)\nearrow 1$ when $T,R\nearrow\infty,$ we obtain
\begin{equation}\label{u=0}
    \Iabt{0}{\infty}{\IRn{|u(t,x)|^p}}\leq 0 \implica u\equiv 0,
\end{equation}
a contradiction.
\vspace{10pt}

\item For the critical case $p=p_c$ in this region $\gamma>\frac{n-2}{n},$ we observe that, from \ref{eq:est5I}, that
\[\Gamma(\alpha)\limsup_{R,T\to\infty}\Iabt{0}{T}{\IRn{\omega_T(t)^\beta |u(t,x)|^{p_c}\varphi_R(x)}}<I_1+I_2-I_3+I_4+I_5,\]
where $I_j, j=1,...,5$ are the five integrals in the right-hand side of inequation \ref{eq:est5I}. Due to $\omega_T,\varphi_R\nearrow 1$ when $T,R\to\infty,$ one obtain
\[\Gamma(\alpha)\limsup_{R,T\to\infty}\Iabt{0}{T}{\IRn{|u(t,x)|^{p_c}}}<I_1+I_2-I_3+I_4+I_5<\infty,\]
hence $u\in L^{p_c}\left([0,\infty)\times \R^n\right).$

Since $\varphi_R$ is constant inside the ball $B_{\frac R2},$ it follows that both $\Delta \varphi_R\equiv0$ and $\Delta^2 \varphi_R\equiv 0$ for $|x|\leq \frac{R}{2}.$ Therefore, integrals $I_1, I_2$ and $I_5$ vanish inside $B_{\frac{R}{2}}.$

On the other hand, outside the ball $B_{R/2}.$ since $u\in L^{p_c}\left([0,\infty)\times \R^n\right),$ we can apply the Riemann-Lebesgue Lemma to see that all $I_1, I_2, I_5$ vanish at infinity.

For $I_3,$ observe that (see Lemma \ref{lema:delTpsi}),
\[\partial_t\psi_T\to 0,\qquad\textrm{as\,\,} T\to\infty.\] Therefore, setting $T=R^{2(1-\theta) }K^{-1},$ for $K\gg1,$
\[\lim_{R\to\infty }I_3=\lim_{R\to\infty }\displaystyle\Iabt{0}{T}{\partial_t\psi_T\IRn{u((-\Delta)^\theta\varphi_R)}} =0.\]

We stress that the idea here was to use the fact that $u\in L^p$ to reestimate and get rid of all the integrals for which the function $g(\eta)$ defined in \eqref{def:g} is attained at the point $\eta=2(1-\theta).$ The estimate for $I_4$ yields the term $(3-\gamma)\eta$ in that definition, and therefore one have some margin at $\eta=2(1-\theta).$ Thus, we have the following estimate:
\begin{align}
\Gamma(\alpha)\limsup_{R,T\to\infty}\Iabt{0}{T}{\IRn{|u(t,x)|^{p_c}}}&\leq \Bigg| \Iabt{0}{\infty}{\partial_{tt}\psi_T\IRn{u\varphi_R}}\Bigg|\notag \\&\leq \varepsilon I(u)+R^{n+2(1-\theta)-2(\alpha+2)(1-\theta)p_c'}K^{-1+(\alpha+2)p_c'},
\end{align}
and so
\[(\Gamma(\alpha)-\varepsilon)\Iabt{0}{\infty}{\IRn{|u(t,x)|^p}}\le C R^{n+2(1-\theta)-2(\alpha+2)(1-\theta)p_c'}.\]

Since $2(1-\theta)(\alpha+2)>g(2(1-\theta)),$ and recalling that
\[p_c'=\dfrac{n+2(1-\theta)}{g(2(1-\theta))},\]
the exponent $n+2(1-\theta)-2(\alpha+2)(1-\theta)p_c'$ is negative, and we have again a contradiction,
\[\Iabt{0}{\infty}{\IRn{|u(t,x)|^p}}\le 0 \implica u\equiv 0.\]

\vspace{10pt}

\item Now, assume that $\gamma\leq\frac{n-2}{n}.$ In this case, the auxiliary function $h(\eta)$ is always non-decreasing. Since
\[\lim_{\eta\to\infty}h(\eta)=\gamma^{-1},\]
it is expected that $\gamma^{-1}$ will replace $p_c$ in this case. Below, we show both cases for $h(\eta);$ the first shows the case where $p_c$ is the maximum of $h(\eta),$ that is, when $\gamma> \frac{n-2}{n};$ the second illustrates the case where $h(\eta)$ does not achieve its maximum, but monotonically tends to its supremum $\gamma^{-1}.$

\begin{minipage}{1.0\textwidth}
\begin{minipage}{0.5\textwidth}
    \includegraphics[scale=0.31]{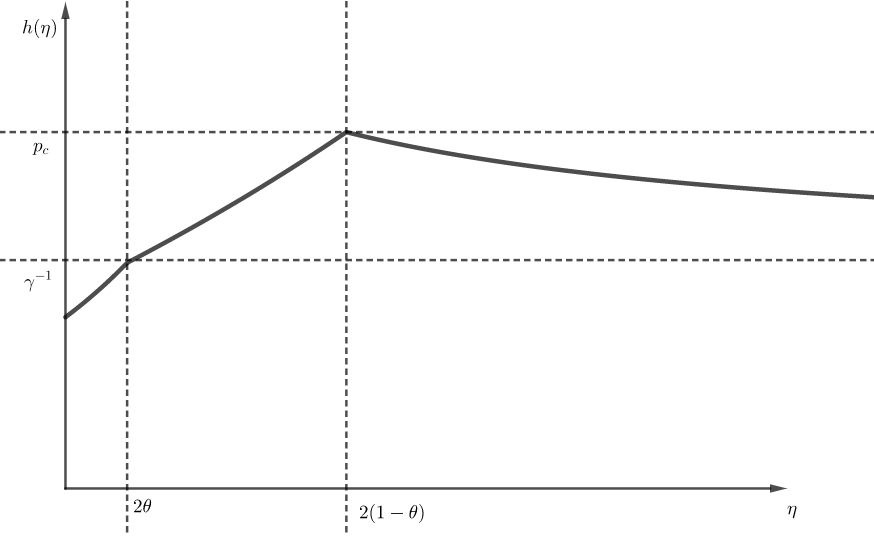}
    \captionof{figure}{$h(\eta),$ $\gamma > \frac{n-2}{n}$}
    \end{minipage}
\begin{minipage}{0.5\textwidth}
    \includegraphics[scale=0.31]{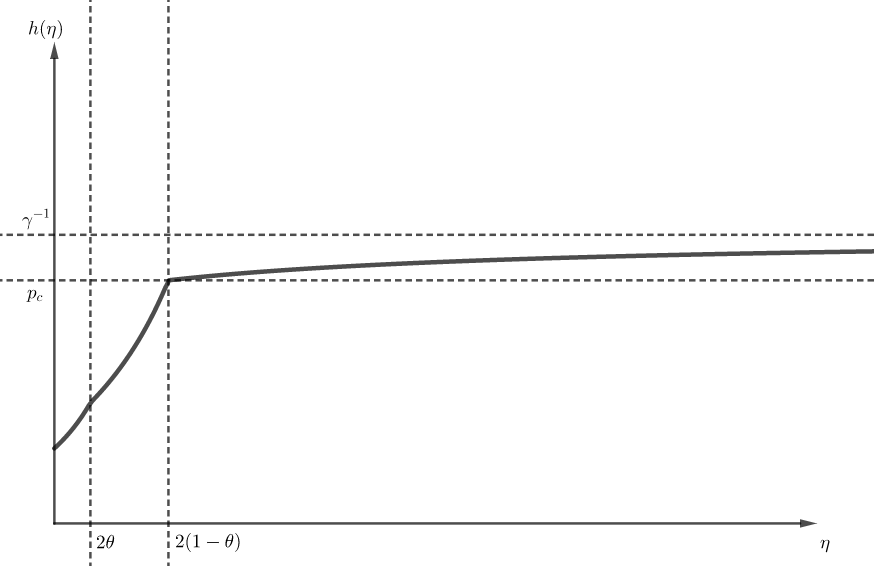}
    \captionof{figure}{$h(\eta),$ $\gamma \leq \frac{n-2}{n}$}
\end{minipage}
\end{minipage}

\bigskip
We set $R=\ln T$ in \eqref{est:nonexistence}, thus obtaining

\small\begin{align}
     \Iabt{0}{T}{\IBR{\omega_T(t)^\beta &|u(t,x)|^p\varphi_R(x)}}\notag\\&< C_{\alpha,\varepsilon}T^{1-\alpha p'}\Big( (\ln T)^{n-2p'}+(\ln T)^{n-4p'}+T^{-p'}(\ln T)^{n-2\theta p'}\notag\\&+T^{-2 p'}(\ln T)^{n}+T^{-2 p'}(\ln T)^{n-2p'}\Big).
\end{align}

\normalsize For $T$ large enough, the logarithmic terms are bounded by $C_\delta T^\delta,$ for any $\delta>0$. Hence,
\begin{align}
     \Iabt{0}{T}{\IBR{\omega_T(t)^\beta |u(t,x)|^p\varphi_R(x)}}&< C_{\alpha,\varepsilon,\delta}T^{1-\alpha p'+\delta}.
\end{align}
If $p<\frac{1}{\gamma},$ then $\alpha p'>1.$ Then, choosing $\delta\in (0,\alpha p'-1),$ we obtain again \eqref{u=0}, which contradicts the assumption that $u$ is nontrivial.

\end{enumerate}

\end{prova}

\end{document}